\documentclass[12pt]{article}
\usepackage{pifont}
\usepackage{bbding}
\usepackage{amssymb}
\usepackage{amsmath}
\usepackage{graphicx}

\usepackage{amsfonts,amssymb,eucal,amsmath}
\title{New Periodic Solutions for Second Order Hamiltonian
Systems with Local Lipschitz Potentials}

\author{\small{\sc    Li Bingyu and Li Fengying and Zhang Shiqing}\\
{\it $\begin{array}{c}
{\rm Mathematical Department,Sichuan University ,Chengdu610064,China}\\
{\rm (Dedicated to the Memory of Professor Shi Shuzhong)}
\end{array}$}}

\date{}

\begin{document}
\maketitle

{\bf Abstract} Firstly,we generalize the classical Palais-Smale-Cerami
condition for $C^1$ functional to the local Lipschitz case,then generalize
 the famous  Benci-Rabinowitz's and Rabinowitz's Saddle Point
Theorems with classical Cerami-Palais-Smale condition to the local Lipschitz functional,
then we apply these Theorems to study the
existence of new periodic solutions  for second
order Hamiltonian  systems with local Lipschitz potentials which are weaker than
Rabinowitz's original conditions .The
key point of our proof is proving Cerami-Palais-Smale
condition for local Lipschitz case,which is difficult since no smooth and symmetry for the potential.\\



\section*{1. Introduction}
\setcounter{section}{1} \setcounter{equation}{0}

In the critical point theory,the compactness condition is a key for proving the existence
  of critical points for some functionals.In 1964,R.Palais and S.Smale [13]introduced the famous $(PS)_c$ condition:

{\bf Definition 1.1}\ \ Let $X$ is a Banach space, $f\in C^1(X,R)$, if $\{x_n\}\subset X$ s.t.

 $$f(x_n)\rightarrow c,$$
 $$f'(x_n)\rightarrow 0,$$

and $\{x_n\}$ has a strongly convergent subsequence,
then we say $f$ satisfies $(PS)_c$ condition.\\

In 1978,Cerami[4] presented a weaker compactness condition than the above
classical $(PS)_c$ condition:

 {\bf Definition 1.2}\ \ Let $X$ be a  Banach
space, $\Phi $ be defined on $X$ is Gateaux-differentiable,
 if the sequence $\{x_n\}\subset X$ such that
 $$\Phi(x_n)\rightarrow c,$$
 $$(1+\|x_n\|)\|\Phi^{'}(x_n)\|\rightarrow 0,$$

then $\{x_n\}$ has a strongly convergent subsequence in $X$.
Then we call $f$ satisfies $(CPS)_c$ condition in $X$.\\

For the functional $f(x)$ in locally Lipschitz functional space $C^{1-0}(X,R)$,Clarke [6] define the generalized gradient $\partial f(x)$ which is the subset of $X^{*}$ defined by
$$\partial f(x)=\{x^{*}\in X^{*} | \langle x^{*} ,v\rangle\leq f^{0}(x,v),\forall v\in X\},$$
where
$$f^{0}(x,v)=\lim_{y\rightarrow x,\lambda\downarrow 0}sup\frac{f(y+\lambda v)-f(y)}{\lambda}.$$

In 1981,K.C.Chang[5] introduced the (PS) condition for locally Lipschitz function:

{\bf Definition 1.3}\ \ Let $X$ is a Banach space, $f\in
C^{1-0}(X,R)$, if $\{x_n\}\subset X$ s.t.$f(x_n)$ is bounded and
 $$min_{x^{*}\in\partial f(x_n)}||x^{*}||\rightarrow 0,$$
and $\{x_n\}$ has a strongly convergent subsequence,
then we say $f$ satisfies $(PSC)$ condition.\\

if $\{x_n\}\subset X$ s.t.$f(x_n)\rightarrow c$  and
 $$min_{x^{*}\in\partial f(x_n)}||x^{*}||\rightarrow 0,$$

and $\{x_n\}$ has a strongly convergent subsequence,
then we say $f$ satisfies $(PSC)_c$ condition.\\

Ekeland [8],Ghoussoub-Preiss[9] used Ekeland's variational principle to prove

{\bf Lemma1.1}Let $X$ be a  Banach space, suppose that $\Phi $ defined on $X$
 is Gateaux-differentiable and lower semi-continuous and bounded from
 below.Then there is a sequence $\{x_n\}$ such that
 $$\Phi(x_n)\rightarrow\inf\Phi$$
 $$(1+\|x_n\|)\|\Phi^{'}(x_n)\|\rightarrow 0.$$

 Motivated by the above Definitions and Lemma,we introduce the following (CPS)-type
 condition for the locally Lipschitz functional:\\

 {\bf Definition 1.4}\ \ Let $X$ is a Banach space, $f\in
C^{1-0}(X,R)$, we say $f$ satisfies $(CPSC)_c$ condition if $\{x_n\}\subset X$ s.t.\\
  $$f(x_n)\rightarrow c,$$
  $$(1+||x_n||)min_{x^{*}\in\partial f(x_n)}||x^{*}||\rightarrow 0,$$
then $\{x_n\}$ has a strongly convergent subsequence.\\

  K.C.Chang[5] and Shi S.Z.[16] use the $(PSC)$ condition for the local Lipschitz functional to generalize the classical Mountain Pass Lemma[2] and general minimax Theorems[12].
  Here we can generalize the classical Benci-Rabinowitz's and Rabinowitz's Saddle Point
Theorems to the local Lipschitz functional cases  with the Cerami-Palais-Smale-Chang-type conditions:

{\bf Theorem1.1}\ \ Let $X$ be a Banach space, $f\in C^{1-0}(X,R)$.
 Let $X=X_1\bigoplus X_2,\rm{dim}X_1<+\infty $,$X_2$ is closed in $X$.
 Let
\begin{eqnarray*}
&&B_a=\{x\in X| \|x\|\leq a\},\\
&&S=\partial B_{\rho}\cap X_2,\rho >0,\\
Q&=&\{x_1+se|(x_1,s)\in X_1\times R^1, \|x_1\|\leq r_1,0\leq s\leq r_2,r_2>\rho\},\\
 &&\partial Q=(B_{r_1}\cap X_1)\cup \partial
\{x_1\bigoplus se, \|x_1\|\leq r_1,0<s\leq r_2\},
\end{eqnarray*}
where $e\in X_2,\|e\|=1$.
If
$$f|_S\geq\alpha ,$$
and
$$f|_{\partial Q}\leq \beta<\alpha.$$

Then $c=\inf\limits_{\phi\in\Gamma}\sup\limits_{x\in Q}f(\phi
(x))\geq\alpha$ ,if $f(q)$ satisfies $(CPSC)_c$
,then $c$ is a critical value for $f$.

{\bf Theorem1.2}\ \ Let $X$ be a Banach space and let $f\in
C^{1-0}(X,R)$, let $X=X_1\bigoplus X_2$ with
$$\rm{dim}X_1<+\infty$$
and
$$\sup\limits_{S^1_R}f<\inf\limits_{X_2}f,$$
where $S^1_R=\{u\in X_1| |u|=R\}$.

Let $B^1_R=\{u\in X_1,|u|\leq R\}, M=\{g\in C(B^1_R,X)|g(s)=s$,
 $s\in S^1_R\}$
$$c=\inf\limits_{g\in M}\max\limits_{s\in B_R^1}(g(s)).$$

Then $c\geq \inf\limits_{X_2}f$, if $f$ satisfies $(CPSC)_c$
condition, then $c$ is a critical value of $f$.\\

In 1978, Rabinowitz [14] firstly used mini-max methods with the classical Palais-Smale condition to
study the periodic solutions for second order Hamiltonian systems
with the super-quadratic potential:

\begin{equation}
\ddot{q}+V^{\prime}(q)=0\label{1.1}
\end{equation}

He proved that

{\bf Theorem 1.3}([14])\ Suppose $V$ satisfies

$(V_1)\ V\in C^1(R^n,R)$

$(V_2)$\ There exist constants $\mu >2, r_0>0$ such that
$$0<\mu V(x)\leq V^{\prime}(x)\cdot x,\ \ \ \ \forall |x|\geq
r_0,$$

$(V_3)\ V(x)\geq 0,\ \ \ \ \forall x\in R^n,$

$(V_4)\ V(x)=o(|x|^2),$\ as\ $|x|\rightarrow 0$.\\

   Then for any $T>0,$ (1.1) has a non-constant $T$-periodic
solution.

In the last 30 years, there were many works for (1.1), we can
refer ([3]-[12],[15,17] etc.), and the references there. In this paper,
we try to generalize the result of Rabinowitz to local Lipschitz potential,
we get the following Theorem:

{\bf Theorem 1.4}\ Suppose $V$ satisfies

$(V1)\ V\in C^{1-0}(R^n,R);$

$(V2)$\ There exist constants ${\mu}_1 >2, \mu_2\in R$ such that
$$ \langle y,x\rangle\geq \mu_1V(x)+\mu_2,\ \ \ \ \forall x\in
R^n ,y\in\partial V(x);$$

$(V3)$\ There are $a_1>0,a_2\in R$ such that
$$V(x)\geq a_1|x|^{\mu_1}+a_2,\ \ \ \ \forall x\in R^n,$$

$(V4)$\ $$ 0\leq V(x)\leq A|x|^2,|x|\rightarrow 0.$$

Then for any $T<(\frac{2}{A})^{1/2}\pi,$ the following system
\begin{equation}
0\in\ddot{q}+\partial V(q)\label{1.2}
\end{equation}

has at least one non-zero $T$-periodic solution.

 For sub-quadratic second order Hamiltonian system,we can get

 {\bf Theorem 1.5}\ Suppose $V$ satisfies

$(V1)\ V\in C^{1-0}(R^n,R);$

$(V2')$\ There exist constants ${\mu}_1<2, \mu_2\in R$ such that
$$ \langle y, x\rangle\leq \mu_1V(x)+\mu_2,\ \ \ \ \forall x\in
R^n ,y\in\partial V(x);$$

$(V3')$\ $$ V(x)\rightarrow +\infty,|x|\rightarrow +\infty;$$

$(V4')$\ $$ V(x)\leq A|x|^2+a.$$

 Then for any
$T<(\frac{2}{A})^{1/2}\pi,$ (1.2) has at least one $T$-periodic solution.

\section*{2. Some Lemmas}
\setcounter{section}{2} \setcounter{equation}{0}

\noindent

     In order to prove Theorem 1.1, we define functional:\\
\begin{equation}
f(q)=\frac{1}{2}\int^T_0|\dot{q}|^2dt-\int^T_0V(q)dt,\ \ \ \
\forall q\in H^1\label{2.1}
\end{equation}

where \begin{equation} H^1=W^{1,2}(R/TZ,R^n).\label{(2.2)}
\end{equation}

{\bf Lemma 2.1([6])}\ \ Let $\widetilde{q}\in H^1$ be such
that $\partial f(\widetilde{q})=0.$\\
Then $\widetilde{q}(t)$ is a  $T$-periodic solution for (1.2).

  {\bf Lemma2.2}(Sobolev-Rellich-Kondrachov, Compact Imbedding
 Theorem [1])\\

$$W^{1,2}(R/TZ,R^n)\subset C(R/TZ,R^n)$$
and the imbedding is compact.

{\bf Lemma 2.3}(Eberlein-Shmulyan [18])\ \ A Banach space $X$ is
reflexive if and only if any bounded sequence in $X$ has a weakly
convergent subsequence.

{\bf Lemma 2.4}([11],[19])\ \ Let $q\in W^{1,2}(R/TZ,R^n)$ and
$q(0)=q(T)=0$

 We have Friedrics-Poincare's inequality:
$$\int^T_0|\dot{q}(t)|^2dt\geq\left(\frac{\pi}{T}\right)^2\int^T_0|q(t)|^2dt.$$\\
Let $q\in W^{1,2}(R/TZ,R^n)$ and $\int^T_0q(t)dt=0,$ then

(i)\
 We have Poincare-Wirtinger's inequality

$$\int^T_0|\dot{q}(t)|^2dt\geq\left(\frac{2\pi}{T}\right)^2\int^T_0|q(t)|^2dt$$

 (ii)\ We have Sobolev's inequality

$$\max_{0\leq t\leq
T}|q(t)|=\|q\|_{\infty}\leq\sqrt{\frac{T}{12}}\left(\int^T_0|\dot{q}(t)|^2dt\right)^{1/2}$$

We define the equivalent norm in $H^1=W^{1,2}(R/TZ,R^n)$
$$\|q\|_{H^1}=\left(\int^T_0|\dot{q}|^2dt\right)^{1/2}+|q(0)|$$

Shi Shuzhong[16] generalized the classical Mini-max Theorems including Benci-Rabinowitz's Generalized
Mountain-Pass Lemma and Rabinowitz's Saddle Point Theorem to the local Lipschitz functionals with Chang's
compactness condition:\\

 {\bf Lemma 2.5}\ \ Let $X$ be a Banach space, $f\in C^{1-0}(X,R)$.
 Let $X=X_1\bigoplus X_2,\rm{dim}X_1<+\infty $,$X_2$ is closed in $X$.Let\\
\begin{eqnarray*}
&&B_a=\{x\in X| \|x\|\leq a\},\\
&&S=\partial B_{\rho}\cap X_2,\rho >0,\\
Q&=&\{x_1+se|(x_1,s)\in X_1\times R^1, \|x_1\|\leq r_1,0\leq s\leq r_2,r_2>\rho\},\\
 &&\partial Q=(B_{r_1}\cap X_1)\cup \partial
\{x_1\bigoplus se, \|x_1\|\leq r_1,0<s\leq r_2\},
\end{eqnarray*}

where $e\in X_2,\|e\|=1$.If
$$f|_S\geq\alpha ,$$
and
$$f|_{\partial Q}\leq \beta<\alpha,$$
Then $c=\inf\limits_{\phi\in\Gamma}\sup\limits_{x\in Q}f(\phi
(x))\geq\alpha$ ,if $f(q)$ satisfies $(PSC)_c$ ,then $c$ is a critical value for $f$.

{\bf Lemma 2.6}\ \ Let $X$ be a Banach space and let $f\in
C^1(X,R)$, let $X=X_1\bigoplus X_2$ with
$$\rm{dim}X_1<+\infty$$
and
$$\sup\limits_{S^1_R}f<\inf\limits_{X_2}f,$$
where $S^1_R=\{u\in X_1| |u|=R\}$.

Let $B^1_R=\{u\in X_1,|u|\leq R\}, M=\{g\in C(B^1_R,X)|g(s)=s$,
 $s\in S^1_R\}$
$$c=\inf\limits_{g\in M}\max\limits_{s\in B_R^1}(g(s))$$
Then $c\geq \inf\limits_{X_2}f$, if $f$ satisfies $(PSC)_c$
condition, then $c$ is a critical value of $f$.\\

  {\bf Lemma 2.7} Let $X$ be a  Banach space, suppose that $F$ defined on $X$
 is local Lipschitz functional and lower semi-continuous and bounded from
 below.Then $\forall\epsilon_n\downarrow 0$, there is a sequence $\{g_n\}$ such that
 $$F(g_n)\rightarrow\inf F,$$
 $$(1+\|g_n\|)F^{0}(g_n,h)|\geq-\epsilon_n\|h\|.$$

{\bf Proof}\ \ Applying Ekeland's variational principle ([7,8]),we can get a sequence $g_n$
such that
$$F(g_n)\leq\inf F+\epsilon_n^2,$$
$$F(g)\geq F(g_n)-\epsilon_n\delta(g,g_n).$$
Let $g=g_n+th,t>0,h\in X$,then we have
$$F(g_n+th)-F(g_n)\geq-\epsilon_n\delta(g_n+th,g_n),$$\\
where $\delta$ is the geodesic distance.
$$F(g_n+th)-F(g_n)\geq-\epsilon_n\int_0^t\frac{||h||ds}{1+||g_n+sh||},$$
then
$$\frac{1}{t}F(g_n+th)-F(g_n)\geq-\epsilon_n\frac{1}{t}\int_0^t\frac{||h||ds}{1+||g_n+sh||},$$
let $t\rightarrow 0$,we have
$$F^{0}(g_n,h)\geq\lim_{t\rightarrow 0}\frac{1}{t}(F(g_n+th)-F(g_n))$$
$$\geq-\epsilon_n\lim_{t\rightarrow 0}\frac{1}{t}\int_0^t\frac{||h||ds}{1+||g_n+sh||}$$
$$=-\epsilon_n||h||(1+||g_n||)^{-1}.$$

\section*{3. The Proof of Theorems 1.1,1.2,1.4 and 1.5}

\setcounter{section}{3} \setcounter{equation}{0}

    {\bf By Lemma 2.7 and similar arguments of Shi Shuzhong [16],we can prove
 Theorem 1.1 and 1.2}.\\

{\bf Lemma 3.1}\ \ If $(V1)-(V3)$ in Theorem 1.4 hold, then $f(q)$
satisfies the $(Cerami-Palais-Smale-Chang)$ condition on $H^1$.

{\bf Proof}\ \ Let $\{q_n\}\subset H^1$ satisfy
\begin{equation}
f(q_n)\rightarrow c,\ \ \ \ (1+||q_n||)min_{x^{*}\in\partial f(q_n)}||x^{*}||\rightarrow 0,\label{3.1}
\end{equation}

Then we claim $\{q_n\}$ is bounded. In fact,by $f(q_n)\rightarrow
c$, we have
\begin{equation}
\frac{1}{2}\|\dot{q}_n\|^2_{L^2}-\int^T_0V(q_n)dt\rightarrow
c\label{3.2}
\end{equation}

By the definition ,we have

$$<\partial f(q_n),q_n>=\|\dot{q}_n\|^2_{L^2}-\int^T_0(<\partial V(q_n),q_n>)dt$$

By $(V2)$,for any $v\in\partial V(q_n)$,we have
\begin{eqnarray}
\|\dot{q}_n\|^2_{L^2}-\int^T_0<v,q_n>dt\leq\|\dot{q}_n\|^2_{L^2}-\int^T_0[\mu_2+\mu_1V(q_n)]dt
\end{eqnarray}

By (3.2) and (3.3), $\forall x^{*}\in\partial f(q_n)$,we have
\begin{eqnarray}
<x^{*},q_n>&\leq&
a\|\dot{q_n}\|^2_{L^2}+C_1+\delta,n\rightarrow+\infty, \label{3.4}
\end{eqnarray}

where $C_1=c\mu_1-T\mu_2+\delta,\delta>0, a=1-\frac{\mu_1}{2}<0.$

By the above inequality (3.4) and (3.1),we have
$\|\dot{q}_n\|_{L^2}\leq M_1$.
 Then we claim $|q_n(0)|$ is also bounded.
Otherwise, there a subsequence, still denoted by $q_n$, s.t.
$|q_n(0)|\rightarrow +\infty$,since $\|\dot{q}_n\|\leq M_1$,then
\begin{eqnarray}
\min_{0\leq t\leq
1}|q_n(t)|&\geq&|q_n(0)|-\|\dot{q}_n\|_2\rightarrow +\infty,
\rm{as}\ n\rightarrow +\infty\label{3.12}
\end{eqnarray}

We notice that
\begin{eqnarray}
\langle\partial f(q_n),q_n\rangle=\int^T_0[|\dot{q}_n|^2dt-\langle\partial V(q_n),q_n\rangle]dt
\end{eqnarray}
\begin{eqnarray}
=2f(q_n)+\int_0^T[2V(q_n)-\langle\partial V(q_n),q_n\rangle]dt \label{3.11}
\end{eqnarray}

By $(V2)-(V3)$,$\forall y\in\partial V(x)$ we have
$$\langle y, x\rangle-2V(x)\geq(\mu_1-2)V+\mu_2\rightarrow +\infty,|x|\rightarrow +\infty$$

By (3.1) and (3.7),we get a contradiction,so $\|q_n\|=\|\dot{q}_n\|_{L^2}+|q_n(0)|$ is bounded.

By the embedding theorem, $\{q_n\}$ has a weakly convergent
subsequence which uniformly converges to $q\in H^1$.

 Furthermore, by $V\in C^{1-0}$ and the $w^{*}-upper$ semi-continuity,
 it's standard step for the rest proof that the weakly convergent subsequence
 is also strongly convergent to $q\in H^1$.

Now we prove {\bf Theorem 1.4.} In Theorem1.1, we take
$$X_1=R^n, X_2=\{q\in W^{1,2}(R/TZ,R^n), \int^T_0q(t)dt=0\}$$
$$S=\left\{q\in
X_2|\left(\int^T_0|\dot{q}|^2dt\right)^{1/2}=\rho>0\right\},$$
$$
\partial Q=\{x_1\in R^n| |x_1|\leq r_1\}\cup$$
$$\left\{q=x_1+se,x_1\in R^n, e\in
X_2,\|e\|=1,s>0,\|q\|=(r_1^2+r_2^2)^{1/2}>\rho\right\}.$$

When $q\in X_2$,by Sobolev's
inequality,$\int^T_0|\dot{q}|^2dt\rightarrow 0$ implies
$max|q(t)|\rightarrow 0$.So when$\int^T_0|\dot{q}|^2dt\rightarrow
0$ , $(V4)$ implies
$$V(q)\leq A|q|^2$$

When $q\in X_2$,we have Poincare-Wirtinger inequality,
so when
$$\rho=[\int^T_0|\dot{q}|^2dt]^{\frac{1}{2}}\rightarrow 0$$

We have
$$f(q)\geq\frac{1}{2}\int^T_0|\dot{q}|^2dt-A\int^T_0|q|^2dt$$
$$\geq[\frac{1}{2}-A(2\pi)^{-2}T^2]\rho^2,$$

On the other hand, if $q\in X_1$,and we take $|x_1|\leq r_1$ very small,
then by $(V_4)$, we have
$$f(q)=-\int^T_0V(q)dt\leq 0,|q|\rightarrow 0.$$

If

$$q\in\left\{q=x_1+se,x_1\in R^n, e\in
X_2,\|e\|=1,s>0,\|q\|=(|x_1|^2+s^2)^{1/2}=R=(r_1^2+r_2^2)^{1/2}>\rho\right\},$$

then by $(V3)$ and Jensen's inequality,we have
$$f(q)=\frac{1}{2}s^2-\int^T_0V(x_1+se)dt$$
$$\leq\frac{1}{2}s^2-\int^T_0(a|x_1+se|^{\mu_1}+b)dt$$
$$\leq\frac{1}{2}s^2-[aT^{1-\frac{\mu_1}{2}}(\int^T_0|x_1+se|^2dt)^{\frac{\mu_1}{2}}+bT]$$
$$=\frac{1}{2}s^2-aT^{1-\frac{\mu_1}{2}}[T|x_1|^2+s^2\int^T_0|e(t)|^2dt]^{\frac{\mu_1}{2}}-bT$$

Notice that we can take $r_2$ large enough,then $(|x_1|^2+s^2)^{1/2}=R=(r_1^2+r_2^2)^{1/2}$ is large enough,then
$|x_1|$ or $s$ must be large,so $T|x_1|^2+s^2\int^T_0|e(t)|^2dt$ must be large since $\int^T_0|e(t)|^2>0$,
so that in such case $f(q)<0.$

The rest of the proof for Theorem 1.4 is obvious.

 {\bf Using Theorem 1.2 and similar methods for proving
Theorem 1.4,we can prove Theorem 1.5,here we omit it}

\section*{Acknowledgements}

 The author Zhang Shiqing sincerely thank  the supports of NSF of China
 and the Grant for the Advisors of Ph.D students.

\end{document}